\documentclass[12pt, makeidx]{amsart}
\usepackage{amsthm}
\usepackage{amsfonts}
\usepackage{amssymb}
\newtheorem{thm}{Theorem}[]

\newtheorem{lem}[thm]{Lemma}

\theoremstyle{definition}

\newtheorem{con}[thm]{Construction}

\newtheorem{notation}{Notation}

\theoremstyle{remark}
\newtheorem{rem}[thm]{Remark}

\numberwithin{thm}{section}
\numberwithin{notation}{section}

\newcommand{\FTWO}{{\rm F}_2}
\newcommand{\Fn}{{\rm F}_n}

\newcommand{\aline}{{\mathcal L}}

\newcommand{\Nat}{\mathbb {N}}

\newcommand{\abs}[1]{\left\vert#1\right\vert}

\newcommand{\wc}{\ (w.c.)}

\begin{document}

\title[On primitives and palindromes]{Palindromic primitives and palindromic bases in the free group of rank two}

\author{Adam Piggott}

\address{School of Mathematics and Applied Statistics, University of Wollongong, NSW 2522, Australia}
\thanks{The author is grateful to Peter Nickolas for a discussion which began this paper.}

\begin{abstract}
The present paper records more details of the relationship between primitive elements and palindromes
in $\FTWO$, the free group of rank two.
We characterise the conjugacy classes of primitive elements which contain palindromes as
those which contain cyclically reduced words of odd length. We identify large palindromic
subwords of certain primitives in conjugacy classes which contain cyclically reduced words
of even length.   We show that under obvious conditions on exponent sums, pairs of
palindromic primitives form palindromic bases for $\FTWO$.  Further, we note
that each cyclically reduced primitive element is either a palindrome, or the concatenation of
two palindromes.
\end{abstract}

\maketitle


\begin{notation}
For each natural number $n \geq 2$, let $\Fn$ denote the non-abelian free group of rank $n$, which we identify with the
set of reduced words in the alphabet $A_n := \{x_1, \dots, x_n\}^\pm$.
For elements $w, v \in \Fn$,  we write $w \equiv_n v$ if $w$ and $v$ are equal words, and $w =_n v$ if $w$ and $v$ are
equal elements of $\Fn$.  We write $wv \wc$ for the concatenation of the words
$w$ and $v$ and $wv$ for the product of $w$ and $v$ in $\Fn$.  We write $\abs{w}_n$ for the word-length of $w$ in $A_n$.
Let $\Psi_n:\Fn \to \Fn$ be the map which reverses each word in $\Fn$.  For convenience, we usually omit the subscript
$n$ from $\equiv_n$, $=_n$, $\Psi_n$, $\abs{\cdot}_n$, and we write $x := x_1$ and $y := y_1$ (so $\FTWO$
is the free group on two generators $x$ and $y$).
\end{notation}

Recall that an element $w \in \Fn$ is said to be a \emph{palindrome}
if $\Psi(w) = w$ (that is, ``$w$ reads the same forwards and
backwards'') and \emph{primitive} if it is an element of some basis for $\Fn$. Much is known about the structure of
primitive elements in $\FTWO$
(see, for example, \cite{WhatDoesABasisLookLike}, \cite{ConstructingPrimitiveElements},
\cite{RemarkOnCountingPrimitiveElements}, \cite{Helling}, \cite{PrimitiveWidth})
and indeed primitive elements in free groups of rank greater than two (see, for example, \cite[pp.162-169]{MKS},
\cite{CountingPrimitiveElements}, \cite{AutomorphicOrbits}, \cite{PrimitiveWidth}).
A newly emerging theme in the study of primitive elements in free groups is the relationship between primitive
elements and palindromes (see \cite{Helling}, \cite{PrimitiveWidth}).  The present paper records
more details of this relationship.

Bardakov, Shpilrain and Tolstykh noted in \cite[p.581]{PrimitiveWidth} that each conjugacy class of primitive elements contains an element
$w$ such that either $xwy^{-1}$ is a palindrome or $x^{-1} w y$ is a palindrome.
It is possible to make a rather more explicit statement concerning the existence of palindromes, and
`large' palindromic subwords, in conjugacy classes of primitive elements.

\begin{thm}\label{PrimitivesAndPalindromes}
Let $p$ be a primitive element in $\FTWO$.  Let $X$ be the exponent sum of $x$ in $p$ and let $Y$ be the exponent
sum of $y$ in $p$. Then:
\begin{enumerate}
\item $X + Y$ is odd if and only if the conjugacy class of $p$ contains exactly one palindrome (a \emph{palindromic primitive});
\item if $X + Y$ is even then the conjugacy class of $p$ does not contain a palindrome but does contain
\begin{enumerate}
\item exactly one element of the form $x^\epsilon w \wc$, and
\item exactly one element of the form $y^\delta v \wc$,
\end{enumerate} where $\epsilon, \delta \in \{\pm 1\}$, the sign of $\epsilon$ (resp. $\delta$)
matches the sign of $X$ (resp. $Y$), and $w, v \in \FTWO$ are palindromes of length $\abs{X} + \abs{Y} -1$.
\end{enumerate}
\end{thm}

Osborne and Zieschang \cite{ConstructingPrimitiveElements}
have recorded an efficient algorithm for writing down a primitive element in $\FTWO$ with a given relatively prime pair of
exponent sums. Theorem \ref{PrimitivesAndPalindromes} is proved by observing the symmetries of a diagrammatic expression of Osborne
and Zieshang's construction.  
The examination of conjugacy classes of primitive elements via the corresponding exponent sums harks back to
Nielsen's work of the early 20th century, which includes the well-known result that conjugacy classes of primitive
elements in $\FTWO$ are in one-to-one correspondence with the set of ordered pairs of integers which are relatively prime,
via the map which takes $w \in \FTWO$ to the pair of exponent sums (see, for example, \cite[p.169]{MKS}, \cite{ConstructingPrimitiveElements}).

Our second theorem demonstrates that, provided the necessary conditions on the exponent sum pairs are satisfied,
pairs of palindromic primitives form bases for $\FTWO$:
\begin{thm}\label{PalindromicBases}
Let $A, B, X$ and $Y$ be integers such that $AY-BX \in \{\pm 1\}$, $A+B$ is odd and $X+Y$ is odd.
The unique palindromic primitive $p$ with exponent sum pair $(A, B)$ and the unique
palindromic primitive $q$ with exponent sum pair $(X, Y)$ form a basis $\{p, q\}$ of $\FTWO$ (a \emph{palindromic basis}).
\end{thm}
\noindent The proof of Theorem \ref{PalindromicBases} also involves an examination of Osborne and
Zieschang's construction.

It is trivial to check that if $w \in \FTWO$ is a product of at most two palindromes, then the image of $w$
under an inner automorphism is also a product of at most two palindromes.  Thus Theorem \ref{PrimitivesAndPalindromes}
supplies another proof of the following:
\begin{lem}[Bardakov, Shpilrain and Tolstykh, Lemma 1.6, p.579 \cite{PrimitiveWidth}]\label{PrimitiveProductOfTwo}
Each primitive element in $\FTWO$ is the product of at most two palindromes.
\end{lem}
\noindent Our third theorem indicates a way in which
Lemma \ref{PrimitiveProductOfTwo} is manifest in the reduced words spelling primitive elements in $\FTWO$:
\begin{thm}\label{ExplicitPalindromes}
For each primitive element $w \in \FTWO$ one of the following holds:
\begin{enumerate}
\item $w$ is a palindromic primitive;
\item $w \equiv pq \wc$ for non-trivial palindromes $p$, $q$;
\item $w \equiv a p a^{-1} \wc$ for a non-trivial palindrome $p$ and a non-trivial word $a \in \FTWO$;
\item $w \equiv a pq a^{-1} \wc$ for non-trivial palindromes $p$, $q$ and a non-trivial word $a \in \FTWO$.
\end{enumerate}
\end{thm}

\noindent Theorem \ref{ExplicitPalindromes} follows immediately from Lemma \ref{PrimitiveProductOfTwo} and
the following result:

\begin{lem}\label{ProductOfTwoCyclicallyReducedLemma}
For each natural number $n \geq 2$, a non-palindromic cyclically reduced element $w \in \Fn$ is a
product of two palindromes in $\Fn$  if and only if
$w \equiv pq \wc$ for palindromes $p, q \in \Fn$.
\end{lem}

\noindent As yet, there is no known algorithm to determine the palindromic (or primitive width) of an element in $\FTWO$ (or more
generally, $F_n$)
\cite[Problems 1 and 2, p.2]{PrimitiveWidth}.  Theorem \ref{ExplicitPalindromes} demonstrates that it is easy to
determine whether or not the palindromic width of an element in $\FTWO$ is zero, one or two.

The structure of the present paper is simple: Theorem \ref{PrimitivesAndPalindromes} is the subject of
Section \ref{PPSection}, Theorem \ref{PalindromicBases} is the subject of Section \ref{PBSection}
and Lemma \ref{ProductOfTwoCyclicallyReducedLemma} is the subject of Section \ref{ProductOfTwoSection}.

\section{Palindromic primitives}\label{PPSection}

Recall the following simple procedure, due to Osborne and Zieshang \cite{ConstructingPrimitiveElements}, for writing down
a primitive element in $\FTWO$ with a given relatively prime pair of exponent sums $X$ and $Y$.

\begin{con}[Osborne and Zieshang, $\S$1.1 of \cite{ConstructingPrimitiveElements}]\label{TheConstruction}
Draw $\abs{X}+\abs{Y}$ equally spaced distinguished points $p_1, p_2, \dots, p_{\abs{X} + \abs{Y}}$
(with indices read around the unit circle in the clockwise direction) on the unit circle in $\mathbb{R}^2$.
Let $l_1$ be $x$ if $X \geq 0$ and $x^{-1}$ if $X < 0$.  Let $l_2$ be $y$ if $Y \geq 0$ and $y^{-1}$ if $Y < 0$.
Label with $l_1$ the points $p_1, p_2, \dots, p_{\abs{X}}$, and label with $l_2$ the remaining
distinguished points.  Let $i$ be an integer such that $1 \leq i \leq \abs{X} + \abs{Y}$ (we call $p_i$
the \emph{first point}).  Let $q_1 := p_i$.
Inductively define $q_j$ for $j = 2, \dots, \abs{X} + \abs{Y}$ as follows: let $q_j$ be the $\abs{X}$-th
distinguished point around the circle from $q_{j-1}$ in the clockwise direction.  For each
$j = 1, \dots, \abs{X} +\abs{Y}$, let $a_j$ be the label on the point $q_j$.  The word
$a_1 a_2 \dots a_{\abs{X} +\abs{Y}}$ is primitive with exponent sum pair $X + Y$.
\end{con}

Note that each of the $\abs{X} + \abs{Y}$ distinct cyclically reduced primitive elements in
the conjugacy class corresponds to a particular choice of first point in Construction \ref{TheConstruction}.

\begin{rem}
It appears to have gone unremarked that a very fast algorithm for
determining whether or not a given cyclically reduced word $w \in \FTWO$ is primitive follows from
Construction \ref{TheConstruction}. We may assume that $\abs{X} + \abs{Y} = \abs{w}$ (since otherwise
$w$ is not primitive \cite{WhatDoesABasisLookLike}).
Place the letters of $w$ around the unit circle, placing
the $j$-th letter so that, travelling around the unit circle in the clockwise direction, the distance from $(1, 0)$
to the $j$-th letter is
$$\frac{2 \pi j \min\{\abs{X}, \abs{Y}\} }{\abs{X} + \abs{Y}} \hbox{ units.}$$
If two letters are placed at the same point, then $w$ is not primitive (since $X$ and $Y$ are not relatively
prime).  Otherwise, $w$ is primitive if and only if the occurrences of $x$ (or $x^{-1}$ if $X < 0$) lie in
consecutive places as read around the circle.
\end{rem}

\begin{proof}[Proof of Theorem \ref{PrimitivesAndPalindromes}]
First consider the case that $X+Y$ is odd.  The result is obvious in the case that $X = 0$ or $Y = 0$, so we may assume
that $X, Y \neq 0$. Consider first the case that $X$ is odd and $Y$ is even.  Consider the diagram of Construction \ref{TheConstruction}
and the line $\aline$ which passes through the Origin and the distinguished point $p_{(\abs{X}+1) / 2}$.
The diagram is symmetric about $\aline$, and $p_{(\abs{X}+1) / 2}$ is the only distinguished point which falls on $\aline$.
Symmetry about $\aline$ ensures that the sequence constructed by reading the label on every $\abs{X}$-th distinguished point
around the circle from $p_{(\abs{X}+1) / 2}$ in the clockwise direction, and the
sequence constructed by reading the label on every $\abs{X}$-th distinguished point around the circle from
$p_{(\abs{X}+1) / 2}$ in the anti-clockwise direction, must be identical.  It follows that
if we choose the first point in Construction \ref{TheConstruction} such that $q_{(\abs{X} + \abs{Y} +1)/2} = p_{(\abs{X}+1) / 2}$,
the result of the construction is a palindrome.  The uniqueness of the
palindrome in the conjugacy class is immediate from the fact that $\aline$ is the unique line of symmetry in the diagram.
The case that $Y$ is odd is handled similarly.

Now consider the case that $\abs{X} + \abs{Y}$ is even.  Since $X$ and $Y$ are relatively prime, both $X$ and $Y$ are odd.
A cyclically reduced element in the conjugacy class of $p$ must have even length.  Palindromes are cyclically reduced.
A palindrome of even length must have even exponent sums, hence there is no palindrome in the conjugacy class of $p$.
Again consider the diagram of Construction \ref{TheConstruction}
and the line $\aline$ which passes through the Origin and the distinguished point $p_{(\abs{X}+1) / 2}$.
This time $\aline$ also passes through the distinguished point $p_{\abs{X}+(\abs{Y} + 1) / 2}$.  Again, $\aline$ is a
line of symmetry in the diagram.  As above, it follows that for each choice of first point in Construction \ref{TheConstruction} such that
$q_{(\abs{X}+\abs{Y}+2)/2} \in \{p_{(\abs{X}+1) / 2}, p_{\abs{X}+(\abs{Y} + 1) / 2}\}$,
the result of the construction is a word of length $\abs{X} + \abs{Y}$
such that the terminal subword of length $\abs{X} + \abs{Y} -1$ is a palindrome.  The fact that there are only two such
words in the conjugacy class is immediate from the fact that $\aline$ is the unique line of symmetry in the diagram.
\end{proof}

\begin{rem}\label{FirstPointForPalindromes}
In case $\abs{X}+\abs{Y}$ is odd and $\abs{X}$ is odd, the choice of first point which gives the unique palindromic primitive is
$p_k$ such that
\begin{equation*}
k \equiv \frac{\abs{X}+1}{2} -\frac{(\abs{X}+\abs{Y}-1)\min\{\abs{X}, \abs{Y}\}}{2} \mod (\abs{X}+\abs{Y}).
\end{equation*}
In case $\abs{X}+\abs{Y}$ is odd and $\abs{Y}$ is odd, the choice of first point which gives the unique palindromic primitive
is $p_k$ such that
\begin{equation*}
k \equiv \abs{X} + \frac{\abs{Y}+1}{2} -\frac{(\abs{X}+\abs{Y}-1)\min\{\abs{X}, \abs{Y}\}}{2} \mod (\abs{X}+\abs{Y}).
\end{equation*}
\end{rem}

\section{Palindromic bases}\label{PBSection}

In this section we prove our result concerning palindromic bases of $\FTWO$.

\begin{proof}[Proof of Theorem \ref{PalindromicBases}]
The Theorem is easily verified in case one of $A, B, X$ or $Y$ is 0, so we may
assume that each is non-zero.  In fact, allowing for the action of automorphisms of $\FTWO$, we may assume
without loss of generality that $0 < A < B$, $0 < X < Y$ and $A+B > X+Y$.  Note that
\begin{equation*}
AY - BX = 1 \Rightarrow A(X+Y) - (A+B)X = 1.    \ \ \ \ \ \ \ \ \ \ \ \ \ \ \ \ \ \ \ \ \ \ Eqn(1)
\end{equation*}
Let $p_1$, $p_2$, $\dots$, $p_{X + Y}$ be the distinguished points when Construction \ref{TheConstruction} is
applied to $X$ and $Y$, and let $w_{X, Y}$ be the result when $p_1$ is chosen as the first point.
Let $r_1$, $r_2$, $\dots$, $r_{A + B}$ be the distinguished points when Construction
\ref{TheConstruction} is applied to $A$ and $B$, and let $w_{A, B}$ be the result when $r_1$ is chosen as the
first point.  By \cite[Theorem 1.2, p.18]{ConstructingPrimitiveElements}, $\{w_{A, B}, w_{X, Y}\}$ is a
basis for $\FTWO$.  By \cite[Theorem 1.3, p.18]{ConstructingPrimitiveElements}, $w_{X, Y}$ is an initial
subword of $w_{A, B}$.

Let $j$ be the integer such that $1 \leq j \leq X + Y$ and $$j \equiv 1 + \frac{(A+B)-(X+Y)}{2} \mod (X + Y).$$
Let $v_{X, Y}$ be the result when $p_{1+jX \mod (X+Y)}$ is chosen as the first letter when Construction \ref{TheConstruction} is
applied to $X$ and $Y$.  The $\frac{X+Y-1}{2}$-th distinguished point visited is $p_k$ such that
$$k \equiv 1 + jX + \frac{(X+Y-1)X}{2} \mod (X+Y).$$  Now
\begin{equation*}
\begin{split}
2k & \equiv 2 + 2jX + (X+Y-1)X \mod (X+Y)\\
 & \equiv 2 + \Bigl(2 + (A+B)-(X+Y)\Bigr)X  \\
 &        \ \ \ \ \ \ \ \ \ \ \ \ \ \ \ + (X+Y-1)X \mod (X+Y) \\
 & \equiv 2 + X + (A+B)X \mod (X+Y) \\
 & \equiv 1 + X  \mod (X+Y), \ \ \ \ \ \ \ \ \ \ \ \ \ \ \ \ \ \ \ \ \ \ \ \ \ \ \ \ \ \ \ \ \ \ \ \ \ \ \ \ \ \ \ \ \ Eqn(2)
\end{split}
\end{equation*}
(the final congruence holds since Eqn(1) implies that $(A+B)X \equiv -\!1 \mod (X+Y)$).

Let $v_{A, B}$ be the result when $q_{1+jA \mod (A+B)}$ is chosen as the first
letter when Construction \ref{TheConstruction} is
applied to $A$ and $B$.  The $\frac{A+B-1}{2}$-th distinguished point visited is $q_\ell$ such that
$$\ell \equiv 1 + jA + \frac{(A+B-1)A}{2} \mod (A+B).$$  Now
\begin{equation*}\label{ResultForl}
\begin{split}
2\ell & \equiv 2 + 2jA + (A+B-1)A \mod (A+B) \\
 & \equiv 2 + \Bigl(2 + (A+B)-(X+Y)\Bigr)A \\
 &     \ \ \ \ \ \ \ \ \ \ \ \ \ \ \ + (A+B-1)A \mod (A+B) \\
 & \equiv 2 + A - (X+Y)A \mod (A+B) \\
 & \equiv 1 + A  \mod (A+B), \ \ \ \ \ \ \ \ \ \ \ \ \ \ \ \ \ \ \ \ \ \ \ \ \ \ \ \ \ \ \ \ \ \ \ \ \ \ \ \ \ \ \ \ \ Eqn(3)
\end{split}
\end{equation*}
(the final congruence holds since Eqn(1) implies that $(X+Y)A \equiv 1 \mod (A+B)$).

The fact that $X + Y$ is odd implies that Eqn(2) has a unique solution for $k$ in the range
$1 \leq k \leq X + Y$.  The fact that $A + B$ is odd implies that Eqn(3) has a unique solution for
$\ell$ in the range $1 \leq \ell \leq A + B$.  If $X$ is odd, it follows (from $AY - BX = 1$) that $A$ is even,
$k = \frac{X+1}{2}$ and $\ell = A + \frac{B+1}{2}$. If $X$ is even, it follows that $A$ is odd, $k = X + \frac{Y+1}{2}$
and $\ell = \frac{A+1}{2}$.  In either case, $p_{k}$ and $q_{\ell}$ are the unique distinguished points that lie
on the line of symmetry in the respective diagrams from Construction \ref{TheConstruction}.  It follows that
$v_{X, Y}$ and $v_{A, B}$ are palindromic primitives and
$$(v_{X, Y}, v_{A, B}) = (c^{-1} w_{X, Y} c, c^{-1} w_{A, B} c),$$
where $c$ is the initial subword of $w_{X, Y}$ (and $w_{A, B}$) of length $j$.  Since $\{w_{X, Y}, w_{A, B}\}$ is a basis
of $\FTWO$, so is $\{v_{X, Y}, v_{A, B}\}$.
\end{proof}

\section{Products of two palindromes}\label{ProductOfTwoSection}

In this section we prove Lemma \ref{ProductOfTwoCyclicallyReducedLemma}, which
combines with Lemma \ref{PrimitiveProductOfTwo} to give Theorem \ref{ExplicitPalindromes}.
We first record the following obvious result:

\begin{lem}\label{EasyCaseLemma}
Let $p \in \Fn$ be a palindrome and $w \in \Fn$ a word:
\begin{enumerate}
\item \label{case1} if $pw \wc$ is a reduced
palindrome and $\abs{p} < \abs{w}$ then $w \equiv q p \wc$ for some
palindrome $q$;
\item \label{case2} if $pw \wc$ is a reduced
palindrome and $\abs{p} = \abs{w}$ then $w \equiv p \wc$;
\item \label{case3} if $wp \wc$ is a reduced
palindrome and $\abs{p} < \abs{w}$ then $w \equiv p q \wc$ for some
palindrome $q$;
\item \label{case4} if $wp \wc$ is a reduced
palindrome and $\abs{p} = \abs{w}$ then $w \equiv p \wc$.
\end{enumerate}
\end{lem}

\begin{lem}\label{PalindromePingPong}
Let $p \in \Fn$ be a palindrome and $w \in F$ a word. If $pw \wc$ is
a reduced palindrome, then one of the following statement holds:
\begin{enumerate}
\item[(1)] $p \equiv r (q r)^m$ and $w \equiv qr  \wc$ for some palindromes $q$, $r$ and some $m \geq
0$;
\item[(2)] $w$ is a palindrome and $p \equiv w^m  \wc$ for some $m \in \Nat$.
\end{enumerate}
If $wp \wc$ is
a reduced palindrome, then one of the following statement holds:
\begin{enumerate}
\item[(3)] $w \equiv rq$ and $p \equiv r (q r)^m  \wc$ for some palindromes $q$, $r$ and some $m \geq
0$;
\item[(4)] $w$ is a palindrome and $p \equiv w^m  \wc$ for some $m \in \Nat$.
\end{enumerate}
\end{lem}

\begin{proof}
Let $a$ be the quotient when
$\abs{p}$ is divided by $\abs{w}$. We prove the lemma using induction on $a$.  By
Lemma \ref{EasyCaseLemma} the result holds in case $a = 0$.  Assume the result holds when $a = k$ for some
non-negative integer $k$.
Consider the case that $a = k+1$.   First we prove the inductive step in case $p w \wc$ is a palindrome.
Now $\abs{p} > \abs{w}$ and $p w \wc$ a reduced palindrome implies
that $p \equiv \Psi(w) p_1  \wc$ for some palindrome $p_1$ such that
$\abs{p_1} < \abs{p}$. Thus we have $p_1$ a reduced palindrome,
$\Psi(w) p_1 \wc$ a reduced palindrome and the quotient when $\abs{p_1}$ is divided by $\abs{w}$ is $k$.  The inductive hypothesis
implies that either
\begin{itemize}
\item $\Psi(w) \equiv rq \wc$, $p_1 \equiv r (q r)^b \wc$ for some palindromes $q$, $r$ and some $b \geq
0$ and hence $p \equiv r (q r)^{b+1} \wc$---Statement (1) holds with $m = b+1$; or
\item $w$ is a palindrome, $p_1 \equiv w^b  \wc$ for some $b > 0$ and hence $p \equiv w^{b+1} \wc$---Statement (2) holds
with $m = b+1$.
\end{itemize}
The inductive step is proved similarly in case $w p \wc$ is a palindrome.
\end{proof}

We are now ready to prove Lemma \ref{ProductOfTwoCyclicallyReducedLemma}:

\begin{proof}[Proof of Lemma \ref{ProductOfTwoCyclicallyReducedLemma}]
We prove only the non-trivial direction of implication. Let $w \in \Fn$ be a non-palindromic product
of non-trivial palindromes $r$ and $s$, that is, $w = rs$.

If the first letter of $r$ and the
last letter of $s$ survive the cancellation between $r$ and $s$, the fact that $w$ is
cyclically reduced implies that the first letter of $r$ and the last
letter of $s$ are not mutually inverse.  It follows that the last
letter of $r$ and the first letter of $s$ are not mutually inverse,
and no cancellation occurs between $r$ and $s$.

If the first letter of $r$ does not survive the cancellation between $r$ and $s$, then $s
\equiv r^{-1} w \wc$.  It follows from Lemma
\ref{PalindromePingPong} that either $w \equiv r^{-n}$ for some $n > 0$, or $r \equiv
(cd)^m c \wc$ and $w \equiv cd \wc$ for some $m \geq 0$ and some palindromes
$c$ and $d$.  The hypothesis that $w$ is not a palindrome implies that
$w \equiv c d \wc$ for non-trivial palindromes $c$ and $d$.  The case that
the last letter of $s$ does not survive the cancellation between $r$ and $s$ is handled similarly.
\end{proof}

\bibliographystyle{amsplain}
\bibliography{PalindromicPrimitiveElementsBib}
\end{document}